\documentclass{amsart}

\usepackage[english]{babel}
\usepackage[utf8x]{inputenc}
\usepackage{amsmath, fullpage}
\usepackage{comment}

\newtheorem{theorem}{Hypothesis}

\title{On the submatrices with the best-bounded inverses}
\author{Richik Sengupta and Mikhail Pautov}

\begin{document}

\begin{abstract} 
The following hypothesis was formulated by Goreinov, Tyrtyshnikov, and Zamarashkin in \cite{goreinov1997theory}. If $U$ is $n\times k$ real matrix with the orthonormal columns $(n>k)$, then there exists a submatrix $Q$ of $U$ of size $k\times k$ such that its smallest singular value is at least $\frac{1}{\sqrt{n}}.$ Although this statement is supported by numerical experiments, the problem remains open for all $1<k<n-1,$ except for the case of $n = 4,\ k=2.$ In this work, we provide a proof for the case $k=2$ and arbitrary $n.$ 
\end{abstract}

\maketitle
\thispagestyle{empty}

\section{Introduction}

This problem was initially formulated in \cite{goreinov1997theory} as the hypothesis about the properties of 
orthonormal $k-$frames. In the case of real numbers, one of the equivalent formulations goes as follows. 

\begin{theorem}
\label{th:main}
    For every $n>k>1$ and arbitrary $n\times k$ real matrix with orthonormal columns, there exists a $k\times k$ submatrix such that the spectral norm of its inverse does not exceed $\sqrt{n}.$ 
\end{theorem}
It is worth mentioning that this problem remains open for all nontrivial cases, except for $n=4, k=2$ (see \cite{nesterenko2023submatrices}). In this paper, we provide a proof for the case $k=2$.

\section{Proof for $k=2$}

\begin{proof}
    Fix $k=2$. Note that the statement for $n=3$ is trivial and the proof for $n=4$ is due to Y. Nesterenko \cite{nesterenko2023submatrices}. Assume that the Hypothesis \ref{th:main} is proven for $n-1$ rows and fix the matrix $A \in \mathbb{R}^{n \times 2}: A^TA=I. $ We want to prove that the statement from the hypothesis holds for the matrix $A.$

\subsection{Case A}

At least one row of the matrix $A$ has a small norm. Specifically, $\exists i: A^2_{i1}+A^2_{i2} \le \frac{1}{n}.$ 
    Without loss of generality, let  $i=1.$  
    Then
    \begin{equation}
    A = \begin{pmatrix}
        A_{11}, & A_{12}\\
       A_{21},  & A_{22} \\
       \dots & \dots \\
       A_{n1}, &A_{n2}
    \end{pmatrix}.
    \end{equation}
    Introduce rotation matrix $P \in \mathbb{R}^{2\times2}$ such that

    \begin{equation}
    B = AP =  \begin{pmatrix}
        b, & 0\\
       B_{21},  & B_{22} \\
       \dots & \dots \\
       B_{n1}, &B_{n2}
    \end{pmatrix}, 
    \text{\ where\ } b^2 = A^2_{11} + A^2_{12} \le \frac{1}{n}.
    \end{equation}
    Note that if $b=0,$ then by removing the first row of $B$, we obtain the matrix of size $(n-1) \times 2$ with orthonormal columns. For such a matrix, according to the induction assumption, there exists a $2\times2$ submatrix with the smallest singular value bounded from below by $\frac{1}{\sqrt{n-1}}\geq \frac{1}{\sqrt{n}}.$ Thus, the statement from Hypothesis \ref{th:main} in Case A trivially holds for $b=0.$ 
    
    Further, we assume that $b^2>0.$

    First,  since singular values of any $2\times2$ submatrix of $A$ are invariant under right multiplication by rotation matrix $P,$ proving the statement from Hypothesis \eqref{th:main} for $B$ is equivalent to proving it for $A.$

    Consider now the submatrix $C$ of the matrix $B,$ where 

    \begin{equation}
        C = \begin{pmatrix}
       B_{21},  & B_{22} \\
       \dots & \dots \\
       B_{n1}, &B_{n2}
    \end{pmatrix}, 
    \end{equation}
where 

\begin{align}
    \begin{cases}
    \sum_{i=2}^n B_{i1}B_{i2} = 0, \\
        \sum_{i=2}^n B^2_{i2} = 1, \\
        \sum_{i=2}^n B^2_{i1} = 1 -b^2.
    \end{cases}
\end{align}
 
Now, if we multiply the first column of $C$ by $t \in \mathbb{R}: t^2\sum_{i=2}^n B^2_{i1} = 1$, or, equivalently, $t^2 = \frac{1}{1-b^2},$ we get the matrix 

\begin{equation}
    \hat C = 
    \begin{pmatrix}
       tB_{21},  & B_{22} \\
       \dots & \dots \\
       tB_{n1}, &B_{n2}
    \end{pmatrix}: \quad \hat C^T \hat C = I.
\end{equation}
Note that $t^2 > 1$ since $b^2>0.$

Now, since $\hat C^T \hat C = I$ and $\hat C$ is $(n-1) \times 2$ matrix, where exists a $2\times 2$ submatrix $\tilde{C}$ of $\hat C$ such that $\sigma^2_2(\tilde{C}) \ge \frac{1}{n-1}.$

Let 
\begin{equation}
    \tilde{C} = \begin{pmatrix}
        tB_{i1} &B_{i2} \\
        tB_{j1} & B_{j2}
    \end{pmatrix} = 
    \begin{pmatrix}
        B_{i1} & B_{i2} \\
        B_{j1} & B_{j2}
    \end{pmatrix} 
    \begin{pmatrix}
        t & 0 \\
        0 & 1
    \end{pmatrix} 
\end{equation}

For $2\times2$ matrices $Y$ and $Z$, we have the inequality for the smallest singular value:
\begin{equation}
    \sigma_2(YZ)\geq\sigma_2(Y)\sigma_2(Z) 
\end{equation}

Thus, 

\[
\sigma_2\begin{pmatrix} B_{i1} & B_{i2} \\ B_{j1} & B_{j2} \end{pmatrix} \geq \sigma_2(\tilde C)\sigma_2\begin{pmatrix} \frac{1}{t} & 0 \\ 0 & 1 \end{pmatrix} \geq \frac{1}{\sqrt{n-1}} \times \frac{1}{t}.
\]

This implies 
\[
\sigma^2_2\begin{pmatrix} B_{i1} & B_{i2} \\ B_{j1} & B_{j2} \end{pmatrix} \geq \frac{1}{n-1}\times \frac{1}{t^2} = \frac{1-b^2}{n-1} \ge \frac{1 - \frac{1}{n}}{n-1} = \frac{1}{n}. 
\]

Thus, there exists a submatrix $\tilde B = \begin{pmatrix} B_{i1} & B_{i2} \\ B_{j1} & B_{j2} \end{pmatrix}$ of $B$ such that $\sigma_2(\tilde B) \ge 1/\sqrt{n},$ and, hence, the statement from Hypothesis  \eqref{th:main} is proven in Case A.

\subsection{Case B} All the norms of the rows of the matrix $A$ are greater than $\frac{1}{\sqrt{n}}:$ $$A^2_{i1}+A^2_{i2}>\frac{1}{n}\  \text{ for all } i\in \overline{1,n}.$$

Let us denote the $i$-th row of the matrix $A$ as $r_i$, i.e., $r_i=(x_i,y_i)$ where $x_i = A_{i1}$ and $y_i=A_{i2}$.

Note that 
\begin{align}
    \sum_{i=1}^n x_i^2 = \sum_{i=1}^n y_i^2 = 1, \quad \sum_{i=1}^n x_iy_i=0. 
\end{align}

We have 

\begin{equation}\label{eq:normsq2}
    \sum_{i=1}^{n} \|r_i\|^2 = \text{Tr}(AA^T) =\text{Tr}(A^T A) = \text{Tr}(I) = 2
\end{equation}

First, we want to prove that there exist $i, j$ such that:

\begin{equation}\label{ineq:main}
(r_i, r_j)^2 \le \left(\|r_i\|^2 - \frac{1}{n}\right)\left(\|r_j\|^2 - \frac{1}{n}\right) 
\end{equation}

\begin{equation}\label{eq:rirj}
    (r_i, r_j)^2 = (x_i x_j + y_i y_j)^2 = \frac{1}{2}(x_i^2 + y_i^2)(x_j^2 + y_j^2) + \frac{1}{2}\big((x_i^2 - y_i^2)(x_j^2 - y_j^2) + 4x_i y_i x_j y_j\big)
\end{equation}
Let us introduce a new set of vectors $w_i \in \mathbb{R}^2$ as $w_i = (x_i^2 - y_i^2, 2x_i y_i)$. 
Notice two important properties of $w_i$:

\begin{equation}\label{eq:wir4}
    \|w_i\|^2 = (x_i^2 - y_i^2)^2 + 4x_i^2 y_i^2 = (x_i^2 + y_i^2)^2 = \|r_i\|^4
\end{equation}

\begin{equation}\label{eq:sumwi0}
    \sum_{i=1}^n w_i = \left(\sum x_i^2 - \sum y_i^2, 2\sum x_i y_i\right) = (1 - 1, 0) = (0,0)
\end{equation}

\vspace{0.5cm}

This allows us to rewrite \eqref{eq:rirj} as:
\[ (r_i, r_j)^2 = \frac{1}{2}\|r_i\|^2 \|r_j\|^2 + \frac{1}{2}(w_i, w_j) \]

Then, after substitution and multiplication by 2,   \eqref{ineq:main} is equivalent to:
\[ \|r_i\|^2 \|r_j\|^2 + (w_i, w_j) \le 2\|r_i\|^2 \|r_j\|^2 - \frac{2}{n}(\|r_i\|^2 + \|r_j\|^2) + \frac{2}{n^2} \]
\[ (w_i, w_j) \le \|r_i\|^2 \|r_j\|^2 - \frac{2}{n}(\|r_i\|^2 + \|r_j\|^2) + \frac{2}{n^2} \]

By adding $\frac{2}{n^2}$ to both sides, we get:

\begin{equation}
\label{eq:wiwj}
    (w_i, w_j) + \frac{2}{n^2} \le \left(\|r_i\|^2 - \frac{2}{n}\right)\left(\|r_j\|^2 - \frac{2}{n}\right) 
\end{equation}

Let us introduce a  set of numbers, $\{z_i\}_{i=1}^n$, such that $z_i = \|r_i\|^2 - \frac{2}{n}$. Note that from \eqref{eq:normsq2}, we get:

\begin{equation}
\label{eq:sumzi0}
    \sum_{i=1}^n z_i = \sum_{i=1}^n \|r_i\|^2 - 2 = 2 - 2 = 0.
\end{equation}

Then \eqref{ineq:main} is equivalent to:

\begin{equation}
\label{eq:wiwj_zizj}
    (w_i, w_j) - z_i z_j + \frac{2}{n^2} \le 0 
\end{equation}

Assume that for all pairs $(i,j)$ the inequality from \eqref{eq:wiwj_zizj} does not hold.  Then for all pairs $(i, j)$:
\begin{equation}
\label{eq:contradiction_wz}
    (w_i, w_j) - z_i z_j + \frac{2}{n^2} > 0.
\end{equation}

Introduce an $n \times n$ matrix $G$ such that $G_{ij} = (w_i, w_j) - z_i z_j$. 

Then $G=WW^T-zz^T,$ where the $i-$th row of $W$ is $w_i$ and $i-$th entry of the vector $z$ is $z_i.$

The matrices $WW^T$ and $zz^T$ are both symmetric and positive semi-definite. Because $W$ has two columns, $\text{rank}(WW^T) \leq 2$, meaning its third largest eigenvalue is bounded by zero:
\begin{equation*}
    \lambda_3(WW^T) \leq 0
\end{equation*}
For any symmetric matrix $X$ and symmetric positive semi-definite matrix $Y$, since $x^T(X-Y)x\leq x^TXx$ for all $x,$ the Courant-Fischer min-max theorem allows to relate  
the $p$-th largest eigenvalues of matrices $X-Y$ and $X$: 
\begin{equation}
    \lambda_p(X - Y) \leq \lambda_p(X) \text{  for all $p \in \overline{1,n}$.}
\end{equation}
Setting $X = WW^T$, $Y = zz^T$, and $p = 3$, we obtain:
\begin{equation*}
    \lambda_3(G) \leq \lambda_3(WW^T) = 0.
\end{equation*}
Because the third largest eigenvalue of $G$ is non-positive, $G$ can have a maximum of two strictly positive eigenvalues.
Consider the trace of $G,$ using \eqref{eq:wir4} and \eqref{eq:normsq2}:

\begin{align*}
\text{Tr}(G) &= \sum_{i=1}^n (\|w_i\|^2 - z_i^2) \\
&= \sum_{i=1}^n \left(\|r_i\|^4 - \left(\|r_i\|^2 - \frac{2}{n}\right)^2\right) \\
&= \sum_{i=1}^n \left(\frac{4}{n}\|r_i\|^2 - \frac{4}{n^2}\right) \\
&= \frac{4}{n}(2) - \frac{4}{n} \\
&= \frac{4}{n}
\end{align*}

Let $\lambda_1$ denote the largest eigenvalue of the symmetric matrix $G.$ Since all its eigenvalues are real, and it has at most $2$ positive eigenvalues,  $\text{Tr}(G)=\frac{4}{n}$ implies that  
\begin{equation}
\label{eq:lambda1greater}
    \lambda_1 \ge \frac{2}{n}.
\end{equation}


Now, define the matrix $M_{ij} = G_{ij} + \frac{2}{n^2}$. Assuming the contradiction from \eqref{eq:contradiction_wz}, all elements of $M$ are positive ($M_{ij} > 0$).

Let $\mathbf{1}\in\mathbb{R}^n$ be the all-ones vector. By linearity, \eqref{eq:sumwi0} and \eqref{eq:sumzi0} yield $G\mathbf{1} = (WW^T-zz^T)\mathbf{1} = \mathbf{0}$. Hence
\[
M\mathbf{1} = G\mathbf{1} + \frac{2}{n^2}\mathbf{E}\mathbf{1} = \frac{2}{n}\mathbf{1},
\]
where $\mathbf{E}$ is the $n\times n$ matrix of ones. Thus, $\frac{2}{n}$ is an eigenvalue of $M$ with positive eigenvector $\mathbf{1}$. \\

By the Perron--Frobenius theorem, eigenspace  associated to  $\frac{2}{n}$ is one-dimensional, and any other eigenvalue $\lambda$ of $M$ satisfies
\begin{equation}\label{Frob-perr}
|\lambda| < \frac{2}{n}.
\end{equation}

Let $\mathbf{1}^\perp = \{x \in \mathbb{R}^n: x^T\mathbf{1} = \mathbf{0}\}.$ For $v\in\mathbf{1}^\perp$ (so $\mathbf{E}v=\mathbf{0}$), we have $Mv = Gv$. Since $G\mathbf{1}=\mathbf{0}$, the eigenvector $v_1$ of $G$ with eigenvalue $\lambda_1\neq0$ lies in $\mathbf{1}^\perp$ by the spectral theorem for symmetric matrices. Then $\lambda_1$ is also an eigenvalue of $M$ with eigenvector $v_1$, so $|\lambda_1| < \frac{2}{n}$ by \eqref{Frob-perr}. This contradicts \eqref{eq:lambda1greater}, hence the assumption \eqref{eq:contradiction_wz} is false.

Therefore, there must exist at least one pair of indices $(i, j)$ such that the \eqref{eq:wiwj_zizj} (and, consequently, \eqref{ineq:main}) holds.

If $i=j$ then after simplification we obtain $\|r_i\|^2 \leq \frac{1}{2n}$ which contradicts the Case B assumption that $\|r_i\|^2 >\frac{1}{n}$ for all $i.$

Therefore, there must exist at least one pair of indices $(i, j)$ such that $i\neq j$ and \eqref{eq:wiwj_zizj} (and, consequently, \eqref{ineq:main}) holds.




Let us consider the corresponding submatrix:

\[
\tilde{A}  = \begin{pmatrix} r_i \\ r_j \end{pmatrix}_{2\times 2}
\]
Its Gram matrix is
\[
\tilde{G} =\tilde{A}\tilde{A}^T= \begin{pmatrix} (r_i,r_i) & (r_i,r_j) \\ (r_j,r_i) & (r_j,r_j) \end{pmatrix}
\]
and the characteristic polynomial is
\[
P_{\tilde{G}}(\lambda)=\left(\|r_i\|^2 - \lambda\right)\left(\|r_j\|^2 - \lambda\right) - (r_i, r_j)^2
\]

 Since for Case B, $\|r_k\|^2>\frac{1}{n}$ for all $k \in \overline{1,n}$ , $Tr(\tilde{G})=\tilde\lambda_1+\tilde\lambda_2 = \|r_i\|^2+\|r_k\|^2 >\frac{2}{n}.$ Here, $\tilde\lambda_1$ and $\tilde\lambda_2 \le \tilde\lambda_1$ are the real eigenvalues of the symmetric matrix $\tilde{G}$. Thus, $\tilde\lambda_1>\frac{1}{n}.$

The characteristic polynomial of $\tilde{G}$ is a quadratic function of $\lambda$ whose leading coefficient is 1 with root(s) at the two eigenvalues . From \eqref{ineq:main} we know, $P_{\tilde{G}}(\frac{1}{n}) \geq 0.$ Thus, $(\tilde\lambda_1 -\frac{1}{n})(\tilde\lambda_2 -\frac{1}{n})\ge 0$.
Thus, $\tilde\lambda_1>\frac{1}{n}$ implies $\tilde\lambda_2\geq \frac{1}{n}.$

 This  implies that for $2 \times 2$ submatrix $\tilde{A} $ of $A,$ we have $\sigma_2(\tilde{A)} = \sqrt{\tilde\lambda_2}\geq\frac{1}{\sqrt{n}}.$ This finalizes the proof for Case B.

\vspace{0.5cm}
 Finally, by induction over $n,$ the statement from Hypothesis \eqref{th:main} is proven for all $n$ and $k=2.$ 
 \end{proof}

\section{Tightness of the bound}
To prove that the bound is tight for arbitrary $n\ge4$, consider the matrix with $n-2$ rows equal to $X = (a,0)$, one row $Y = (b,c)$, and another row $Z = (b,-c)$. Let 
$$a = \sqrt{\frac{n-1}{n(n-2)}}, \quad b = \frac{1}{\sqrt{2n}}, \quad \text{and} \quad c = \frac{1}{\sqrt{2}}.$$
Specifically, such a matrix can be of the form:

\[
\begin{pmatrix}
a & 0 \\
a & 0 \\
\vdots & \vdots \\
a & 0 \\
b & c \\
b & -c
\end{pmatrix}
\]
Let $\mathcal{A}$ denote the set of all $2 \times 2$ submatrices formed by the rows of this matrix, and let $\sigma_2(\tilde{A})$ denote the smallest singular value of a given submatrix $\tilde{A}$. Direct evaluation demonstrates that:
$$\max_{\tilde{A} \in \mathcal{A}} \sigma_2(\tilde{A}) = \frac{1}{\sqrt{n}}$$
Furthermore, the bound is attained on the following submatrices:
$$\left\{ \begin{pmatrix} X \\ Y \end{pmatrix}_{2\times 2}, \begin{pmatrix} X \\ Z \end{pmatrix}_{2\times 2}, \begin{pmatrix} Y \\ Z \end{pmatrix}_{2\times 2} \right\} \subseteq \operatorname*{arg\,max}_{\tilde{A} \in \mathcal{A}} \sigma_2(\tilde{A}).$$
Thus, the bound is tight in the sense that  $\frac{1}{\sqrt{n}}$ can not be replaced by a larger lower bound that holds uniformly for all real matrices $A \in \mathbb{R}^{n\times 2}: A^TA=I.$

 \section{Conclusion}
 The proven result positively resolves the hypothesis from equation (2.6) for $r=2$ in \cite{goreinov1997theory} and provides an explicit upper bound on errors of pseudoskeleton approximations, which are formulated as Theorems 3.1 and 3.2 in \cite{goreinov1997theory}. More examples of matrices for which the lower bound is  attained are discussed in \cite{nesterenko2024submatrices}.

\bibliographystyle{plain} 
\bibliography{main}

\end{document}